\documentclass{amsart}
\usepackage{amssymb}

\usepackage{graphicx}
\usepackage{amscd}

\theoremstyle{plain}

\newtheorem{lemma}{Lemma}

\numberwithin{equation}{section}

\input{tcilatex}

\begin{document}
\title[Complete minimal hypersurfaces ]{Complete minimal hypersurfaces of $\mathbb{S}^{4}$ with zero Gauss-Kronecker
curvature}
\author{T. Hasanis}
\address[]{Department of Mathematics, University of Ioannina, 45110, Ioannina, Greece}
\email{thasanis@cc.uoi.gr, me00499@cc.uoi.gr, tvlachos@cc.uoi.gr}
\author{A. Savas-Halilaj}
\author{T. Vlachos}

\begin{abstract}
We investigate the structure of 3-dimensional complete minimal hypersurfaces
in the unit sphere with Gauss-Kronecker curvature identically zero.
\end{abstract}

\maketitle

\section{Introduction}

In order to consider rigidity problems, Dajczer and Gromoll in \cite{DG}
studied the Gauss map of hypersurfaces, which in general is not invertible.
The starting point is the observation that whenever the relative nullity is
constant, then one has a representation of the hypersurface by the inverse
of the Gauss map on the normal bundle of its image, which they called ''%
\textit{Gauss parametrization}''. In particular, they proved that if $%
g:V\rightarrow \mathbb{S}^{n+1}$ is a minimal immersion, where $V$ is a
2-dimensional manifold,\ and $\mathcal{N}$ is the unit normal bundle of $g$,
then the ''\textit{polar map'' }$\Psi \left( p,w\right) =w$, $\left(
p,w\right) \in \mathcal{N}$, defines a minimal hypersurface which at its
regular points has relative nullity $n-2$. Conversely, it is shown that,
locally, every minimal hypersurface with constant relative nullity $n-2$ has
such a representation.

Consider a minimal immersion $g:V\rightarrow \mathbb{S}^{4}$, where $V$ is a
2-dimensional manifold, and denote by $\mathcal{K}$, $\mathcal{K}_{n}$ its
Gaussian and normal curvature, respectively. It is proved in \cite{DG} that
if $g$ has nowhere vanishing\ normal curvature, then the polar map of $g$ is
everywhere regular and provides a minimal hypersurface in $\mathbb{S}^{4}$
with Gauss-Kronecker curvature identically zero. In \cite{AB} de Almeida and
Brito classified the compact minimal hypersurfaces $M^{3}$ in $\mathbb{S}%
^{4} $ with Gauss-Kronecker curvature identically zero under the assumption
that the second fundamental form of $M^{3}$ nowhere vanishes. More
precisely, they proved that such a hypersurface is the boundary of a tube,
of radius $\frac{\pi }{2}$, of a minimal immersion $g:V\rightarrow \mathbb{S}%
^{4}$ with non-vanishing second fundamental form in any direction. Later,
Ramanathan in \cite{R} removed the assumption on the second fundamental form
and classified the compact minimal hypersurfaces in $\mathbb{S}^{4}$ with
Gauss-Kronecker identically zero. In particular,\ he proved that such
hypersurfaces are produced by applying the above construction to appropriate
branched minimal surfaces in $\mathbb{S}^{4},$ unless they are totally
geodesic.

A minimal immersion $g:V\rightarrow \mathbb{S}^{4}$ is called \textit{%
superminimal }if it satisfies the relation\textit{\ }$\left( 1-\mathcal{K}%
\right) ^{2}-\mathcal{K}_{n}^{2}=0.$ The purpose of the present paper is to
consider the structure of complete minimal immersions $f:M^{3}\rightarrow 
\mathbb{S}^{4}$ with Gauss-Kronecker curvature identically zero under some
assumptions on the second fundamental form. More precisely, we show the
following

\noindent \textbf{Theorem.} \textit{Let }$f:M^{3}\rightarrow \mathbb{S}^{4}$%
\textit{\ be a minimal isometric immersion with Gauss-Kronecker curvature
identically zero, where }$M^{3}$\textit{\ is a 3-dimensional complete
Riemannian manifold. If the square }$S$\textit{\ of the length of the second
fundamental form is nowhere zero and bounded from above, then }$f\left(
M^{3}\right) $\textit{\ is the image of the polar map associated\ with a
superminimal immersion }$g:V\rightarrow \mathbb{S}^{4}$\textit{\ with
positive normal curvature. Moreover, if }$S$\textit{\ is bounded away from
zero, then }$V$\textit{\ is diffeomorphic to the sphere }$\mathbb{S}^{2}$%
\textit{\ or to the projective plane }$\mathbb{R}\mathbb{P}^{2}$\textit{and }%
$f\left( M^{3}\right) $\textit{\ is compact.}

\section{Preliminaries}

Consider an isometric and minimal immersion $g:V\rightarrow \mathbb{S}^{4},$
where $V$ is an oriented, 2-dimensional Riemannian manifold. Denote by $%
\left( v_{1},v_{2};v_{3},v_{4}\right) $ a local adapted orthonormal frame
field along $g$ such that $\left( v_{1},v_{2}\right) $ is an oriented
orthonormal frame field in the tangent bundle and $\left( v_{3},v_{4}\right) 
$ is an oriented frame field in the normal bundle of $g$. Denote by $A_{3}$, 
$A_{4}$ the shape operators corresponding to the directions $v_{3}$ and $%
v_{4}$. The Gaussian curvature $\mathcal{K}$ of the induced metric $%
\left\langle \ ,\ \right\rangle $ and the normal curvature $\mathcal{K}_{n}$
of the immersion $g$ are given by 
\begin{equation*}
\mathcal{K}=1+\det A_{3}+\det A_{4},\text{ }\mathcal{K}_{n}=-\left\langle %
\left[ A_{3},A_{4}\right] v_{1},v_{2}\right\rangle ,
\end{equation*}
and do not depend on the chosen frames. Then, $\left( 1-\mathcal{K}\right)
^{2}-\mathcal{K}_{n}^{2}\geq 0$. The isometric immersion $g$ is called 
\textit{superminimal} if it\ satisfies in addition $\left( 1-\mathcal{K}%
\right) ^{2}-\mathcal{K}_{n}^{2}=0$. In this case, we can choose the adapted
orthonormal frame field $\left( v_{1},v_{2};v_{3},v_{4}\right) $ such that 
\begin{equation*}
A_{3}\sim \left( 
\begin{array}{ll}
\mu  & \ \ 0 \\ 
0 & -\mu 
\end{array}
\right) ,\text{ \ }A_{4}\sim \left( 
\begin{array}{ll}
0 & \mu  \\ 
\mu  & 0
\end{array}
\right) .
\end{equation*}
It is obvious that the square of the length of the second fundamental form
of $g$ is equal to $2\left( 1-\mathcal{K}\right) $. Tribuzy and Guadalupe
in\ \cite{TG} proved that if $g$ is superminimal then the Gaussian curvature
satisfies the following differential equation 
\begin{equation}
\Delta \log \left( 1-\mathcal{K}\right) =2\left( 3\mathcal{K}-1\right) ,
\end{equation}
away from points where $\mathcal{K}=1$, where $\Delta $ stands for the
Laplacian operator on $V$. Using the Penrose twistor fibration of complex
projective 3-space over $\mathbb{S}^{4}$, Bryant in \cite{B} was able to
construct superminimal immersions of $\mathbb{S}^{2}$ into $\mathbb{S}^{4}$,
even with nowhere vanishing normal curvature.

We give now a brief exposition of a method developed by Dajczer and Gromoll
in \cite{DG} (see also \cite{R}, Theorem 1.1) of constructing minimal
hypersurfaces in $\mathbb{S}^{4}$ with Gauss-Kronecker curvature identically
zero. Let $g:V\rightarrow \mathbb{S}^{4}$ be an isometric minimal immersion,
where $V$ is a 2-dimensional Riemannian manifold, and 
\begin{equation*}
\mathcal{N}=\left\{ \left( x,w\right) \in V\times \mathbb{R}^{5}:\left|
w\right| =1,\text{ }w\perp \mathbb{R\cdot }g+dg\left( T_{x}V\right) \right\}
.
\end{equation*}
be its unit normal bundle. Denote the projection to the first factor by $\pi
:\mathcal{N}\rightarrow V$. The projection to the second factor $\Psi :%
\mathcal{N}\rightarrow \mathbb{S}^{4}$, $\Psi \left( x,w\right) =w$, is
called the ''\textit{polar map''} associated with $g$. Choose an adapted
orthonormal frame field $\left( v_{1},v_{2};v_{3},v_{4}\right) $ on an open
set $U\subset V$ such that the shape operators $A_{3}$, $A_{4}$ of $g$
corresponding to $v_{3}$, $v_{4}$ are represented by 
\begin{equation*}
A_{3}\sim \left( 
\begin{array}{cc}
a & \text{ \ }b \\ 
b & -a
\end{array}
\right) \text{, \ }A_{4}\sim \left( 
\begin{array}{cc}
c & \text{ \ }0 \\ 
0 & -c
\end{array}
\right) \text{,}
\end{equation*}
and parametrize $\pi ^{-1}\left( U\right) $ by $U\times \mathbb{S}^{1}$ via
the map $\left( x,t\right) \rightarrow \left( x,\cos tv_{3}\left( x\right)
+\sin tv_{4}\left( x\right) \right) $. Then, $\Psi \left( x,t\right) =\cos
tv_{3}\left( x\right) +\sin tv_{4}\left( x\right) $. Hence, 
\begin{eqnarray*}
d\Psi \left( \partial /\partial t\right)  &=&-\sin tv_{3}+\cos tv_{4}, \\
d\Psi \left( v_{1}\right)  &=&-\left( a\cos t+c\sin t\right) dg\left(
v_{1}\right) -\left( b\cos t\right) dg\left( v_{2}\right)  \\
&&+\omega _{34}\left( e_{1}\right) \left( -\sin tv_{3}+\cos tv_{4}\right) ,
\\
d\Psi \left( v_{2}\right)  &=&-\left( b\cos t\right) dg\left( v_{1}\right)
+\left( a\cos t+c\sin t\right) dg\left( v_{2}\right)  \\
&&+\omega _{34}\left( e_{2}\right) \left( -\sin tv_{3}+\cos tv_{4}\right) ,
\end{eqnarray*}
where $\omega _{34}$ is the connection form of the normal bundle of $g$.
From the above relations, it follows that $\Psi $ is regular at points $%
\left( x,t\right) $, for all$\ t\in \mathbb{S}^{1}$, if and only if $\left(
a\cos t+c\sin t\right) ^{2}+\left( b\cos t\right) ^{2}\neq 0$ which is
equivalent to $\mathcal{K}_{n}\neq 0$. Obviously, $\xi \left( x,t\right)
=g\left( x\right) $ defines a unit normal vector field along $\Psi $. Using
the Weingarten formulas one get that $\Psi $, where it is regular, has
principal curvatures 
\begin{equation*}
k_{1}=-k_{3}=\frac{1}{\sqrt{\left( a\cos t+c\sin t\right) ^{2}+\left( b\cos
t\right) ^{2}}},\text{ }k_{2}=0\text{.}
\end{equation*}

Let $M^{3}$ be a 3-dimensional, oriented Riemannian manifold and $%
f:M^{3}\rightarrow \mathbb{S}^{4}$ an isometric minimal immersion into the
unit sphere. Denote by $\xi $ a unit normal vector field along $f$, by $A$
the shape operator associated with $\xi $ and by $k_{1}\geq k_{2}\geq k_{3}$
the principal curvatures. The Gauss-Kronecker curvature $K$ and the square $%
S $ of the length of the second fundamental form are given, respectively, by 
\begin{equation*}
K=k_{1}k_{2}k_{3},\ \ S=k_{1}^{2}+k_{2}^{2}+k_{3}^{2}.
\end{equation*}
Assume now that $K=0$ and that the second fundamental form is nowhere zero.
Then the principal curvatures satisfy the relation $k_{1}=:\lambda
>k_{2}=0>k_{3}=-\lambda $, where $\lambda $ is a smooth positive function on 
$M^{3}$. We can choose locally an orthonormal frame field $\left(
e_{1},e_{2},e_{3}\right) $ of principal directions corresponding to $\lambda
,0,-\lambda $. Let $\left( \omega _{1},\omega _{2},\omega _{3}\right) $ and $%
\omega _{ij}$, $i,j\in \{1,2,3\}$, be the corresponding dual and the\
connection forms, respectively. Throughout this paper we make the following
convection on the ranges of indices 
\begin{equation*}
1\leq i,j,k,\ldots \leq 3
\end{equation*}
and adopt the method of moving frames. The structural equations are 
\begin{eqnarray*}
d\omega _{i} &=&\sum_{j}\omega _{ij}\wedge \omega _{j},\quad \omega
_{ij}+\omega _{ji}=0, \\
d\omega _{ij} &=&\sum_{l}\omega _{il}\wedge \omega _{lj}-\left(
1+k_{i}k_{j}\right) \omega _{i}\wedge \omega _{j}.
\end{eqnarray*}
Consider the functions 
\begin{equation*}
u:=\omega _{12}\left( e_{3}\right) ,\quad v:=e_{2}\left( \log \lambda
\right) ,
\end{equation*}
which will play a crucial role in the proof of our result. From the
structural equations, and the Codazzi equations, 
\begin{equation*}
\begin{array}{l}
e_{i}\left( k_{j}\right) =\left( k_{i}-k_{j}\right) \omega _{ij}\left(
e_{j}\right) ,\ i\neq j,\medskip \\ 
\left( k_{1}-k_{2}\right) \omega _{12}\left( e_{3}\right) =\left(
k_{2}-k_{3}\right) \omega _{23}\left( e_{1}\right) =\left(
k_{1}-k_{3}\right) \omega _{13}\left( e_{2}\right) ,
\end{array}
\end{equation*}
we, easily, get 
\begin{equation}
\begin{array}{lll}
\omega _{12}\left( e_{1}\right) =v, & \omega _{13}\left( e_{1}\right) =\frac{%
1}{2}e_{3}\left( \log \lambda \right) , & \omega _{23}\left( e_{1}\right)
=u,\medskip \\ 
\omega _{12}\left( e_{2}\right) =0, & \omega _{13}\left( e_{2}\right) =\frac{%
1}{2}u, & \omega _{23}\left( e_{2}\right) =0,\medskip \\ 
\omega _{12}\left( e_{3}\right) =u, & \omega _{13}\left( e_{3}\right) =-%
\frac{1}{2}e_{1}\left( \log \lambda \right) , & \omega _{23}\left(
e_{3}\right) =-v,
\end{array}
\end{equation}
and 
\begin{equation}
\begin{array}{l}
e_{2}\left( v\right) =v^{2}-u^{2}+1,\ e_{1}\left( u\right) =e_{3}\left(
v\right) ,\ e_{2}\left( u\right) =2uv,\ e_{3}\left( u\right) =-e_{1}\left(
v\right) .
\end{array}
\end{equation}
Observe that the integral curves of $e_{2}$ are geodesics in $M^{3}$ and
their images under $f$ are geodesics in $\mathbb{S}^{4}$. Furthermore, from
the\ above equations, for the Lie bracket we get 
\begin{equation}
\lbrack e_{1},e_{3}]=-\frac{1}{2}e_{3}\left( \log \lambda \right)
e_{1}-2ue_{2}+\frac{1}{2}e_{1}\left( \log \lambda \right) e_{3}.
\end{equation}

\section{Proof Of The Theorem}

Our main result follows from a sequence of lemmas which are themselves of
independent interest.

\begin{lemma}
Under the notation introduced in section 2 we have:$\medskip $

\noindent {}$(i)$ The function $\dfrac{u}{\lambda ^{2}}$ is constant along
the integral curves of $e_{2}$.

\noindent {}$(ii)$ The functions $u$ and $v$ are harmonic.
\end{lemma}

\begin{proof}
$(i)$ By making use of $\left( 2.2\right) $ and $\left( 2.3\right) $ we,
immediately, obtain $e_{2}\left( \frac{u}{\lambda ^{2}}\right) =0$.

\noindent $(ii)$ From the definition of Laplacian we have 
\begin{eqnarray*}
\Delta v &=&e_{1}e_{1}\left( v\right) +e_{2}e_{2}\left( v\right)
+e_{3}e_{3}\left( v\right) -\left( \omega _{21}\left( e_{2}\right) +\omega
_{31}\left( e_{3}\right) \right) e_{1}\left( v\right) \smallskip  \\
&&-\left( \omega _{12}\left( e_{1}\right) +\omega _{32}\left( e_{3}\right)
\right) e_{2}\left( v\right) -\left( \omega _{13}\left( e_{1}\right) +\omega
_{23}\left( e_{2}\right) \right) e_{3}\left( v\right) 
\end{eqnarray*}
or, taking into account $\left( 2.2\right) ,$%
\begin{eqnarray}
\Delta v &=&e_{1}e_{1}\left( v\right) +e_{2}e_{2}\left( v\right)
+e_{3}e_{3}\left( v\right) -\frac{1}{2}e_{1}\left( \log \lambda \right)
e_{1}\left( v\right)  \\
&&-2ve_{2}\left( v\right) -\frac{1}{2}e_{3}\left( \log \lambda \right)
e_{3}\left( v\right) .  \notag
\end{eqnarray}
We also have from $\left( 2.3\right) $ 
\begin{eqnarray}
e_{1}e_{1}\left( v\right)  &=&-e_{1}e_{3}\left( u\right) ,\ e_{3}e_{3}\left(
v\right) =e_{3}e_{1}\left( u\right) ,\smallskip  \\
e_{2}e_{2}\left( v\right)  &=&2ve_{2}\left( v\right) -2ue_{2}\left( u\right)
=2v^{3}-6vu^{2}+2v.
\end{eqnarray}
Inserting $\left( 3.2\right) $ and $\left( 3.3\right) $ into $\left(
3.1\right) $ and using $\left( 2.4\right) ,$ we obtain 
\begin{eqnarray*}
\Delta v &=&-e_{1}e_{3}\left( u\right) +e_{3}e_{1}\left( u\right)
+2v^{3}-6vu^{2}+2v \\
&&-\frac{1}{2}e_{1}\left( \log \lambda \right) e_{1}\left( v\right)
-2ve_{2}\left( v\right) -\frac{1}{2}e_{3}\left( \log \lambda \right)
e_{3}\left( v\right)  \\
&=&\frac{1}{2}e_{3}\left( \log \lambda \right) e_{1}\left( u\right)
+2ue_{2}\left( u\right) -\frac{1}{2}e_{1}\left( \log \lambda \right)
e_{3}\left( u\right)  \\
&&+2v^{3}-6vu^{2}+2v \\
&&-\frac{1}{2}e_{1}\left( \log \lambda \right) e_{1}\left( v\right)
-2ve_{2}\left( v\right) -\frac{1}{2}e_{3}\left( \log \lambda \right)
e_{3}\left( v\right) .
\end{eqnarray*}
Appealing to $\left( 2.3\right) $, we readily see that $v$ is harmonic.

In a similar way, we verify that $\Delta u=0$.
\end{proof}

\begin{lemma}
Let $M^{3}$ be a 3-dimensional, oriented, complete Riemannian manifold and $%
f:M^{3}\rightarrow \mathbb{S}^{4}$ a minimal isometric immersion with
Gauss-Kronecker curvature identically zero and nowhere vanishing second
fundamental form. Then the function $u$ is nowhere zero.
\end{lemma}

\begin{proof}
Assume in the contrary that there exists a point $x_{0}\in M^{3}$ such that$%
\ u\left( x_{0}\right) =0$. Let $\gamma \left( s\right) $, $s\in \mathbb{R}$%
, be the maximal integral curve of $e_{2}$ emanating from the point $x_{0},$
where $s$ is its arclength. Because of Lemma $1$ $\left( i\right) ,$ the
function $u\left( s\right) :=u\left( \gamma \left( s\right) \right) $ must
be everywhere zero. Restricting the first equation of $\left( 2.3\right) $
along $\gamma \left( s\right) ,$ we obtain the differential equation $%
v^{\prime }\left( s\right) =v^{2}\left( s\right) +1,$ where $v\left(
s\right) :=v\left( \gamma \left( s\right) \right) $ is an entire function.
This is a contradiction because this equation cannot admit entire solutions.
\end{proof}

\begin{lemma}
Let $M^{3}$ be a 3-dimensional, oriented, complete Riemannian manifold and $%
f:M^{3}\rightarrow \mathbb{S}^{4}$ a minimal isometric immersion with
Gauss-Kronecker curvature identically zero and nowhere vanishing second
fundamental form. \textit{Let }$\xi :M^{3}\rightarrow \mathbb{S}^{4}$ be the
Gauss map; then there exists a 2-dimensional differentiable manifold $V$, a
submersion $\pi :M^{3}\rightarrow V$ and a minimal immersion $\widetilde{\xi 
}:V\rightarrow \mathbb{S}^{4}$, with nowhere vanishing normal curvature,
such that $\pi \circ \widetilde{\xi }=\xi .$
\end{lemma}

\begin{proof}
Consider the quotient space $V$ of leaves of $e_{2}$ with quotient map $\pi
:M^{3}\rightarrow V$. Since $M^{3}$ is complete, the integral curves of $%
e_{2}$ are complete geodesics and their images through $f$ are great circles
of $\mathbb{S}^{4}$. These facts ensure $\left( \text{see \cite{P}}\right) $
that $V$ can be equipped with a structure of a 2-dimensional differentiable
manifold which makes $\pi $ a submersion. The Gauss map $\xi $ induces a
smooth map $\widetilde{\xi }:V\rightarrow \mathbb{S}^{4}$ so that $%
\widetilde{\xi }\circ \pi =\xi $. Consider, now, a smooth transversal $S$ to
the leaves of $e_{2}$ through a point $x\in M^{3}$ such that $e_{1}\left|
_{x}\right. $, $e_{3}\left| _{x}\right. $ span $T_{x}S$. Because $\pi $ is
submersion,$\ \left( d\pi \left( e_{1}\left| _{x}\right. \right) ,d\pi
\left( e_{3}\left| _{x}\right. \right) \right) $ constitute a base of $%
T_{\pi \left( x\right) }V$. Note that 
\begin{equation*}
d\widetilde{\xi }\left( d\pi \left( e_{1}\left| _{x}\right. \right) \right)
=-\lambda \left( x\right) df\left( e_{1}\left| _{x}\right. \right) ,\text{ }d%
\widetilde{\xi }\left( d\pi \left( e_{3}\left| _{x}\right. \right) \right)
=\lambda \left( x\right) df\left( e_{3}\left| _{x}\right. \right) .
\end{equation*}
Thus $\widetilde{\xi }$ is an immersion and $X_{1}:=d\pi \left( \frac{1}{%
\lambda }e_{1}\left| _{x}\right. \right) ,$ $X_{2}:=d\pi \left( \frac{1}{%
\lambda }e_{3}\left| _{x}\right. \right) $ are orthonormal at $\pi \left(
x\right) $ with respect to the metric induced by $\widetilde{\xi }$. Let $%
N_{1},$ $N_{2}$ be an orthonormal frame in the normal bundle of $\widetilde{%
\xi }$ such that $N_{1}\circ \pi \left| _{S}\right. =f\left| _{S}\right. $, $%
N_{2}\circ \pi \left| _{S}\right. =df\left( e_{2}\right) \left| _{S}\right. $%
. Observe that 
\begin{eqnarray*}
dN_{1}\left( X_{1}\right) &=&\frac{1}{\lambda \left( x\right) }df\left(
e_{1}\left| _{x}\right. \right) ,\text{ }dN_{1}\left( X_{2}\right) =\frac{1}{%
\lambda \left( x\right) }df\left( e_{3}\left| _{x}\right. \right) \medskip ,
\\
dN_{2}\left( X_{1}\right) &=&-\frac{v\left( x\right) }{\lambda \left(
x\right) }df\left( e_{1}\left| _{x}\right. \right) +\frac{u\left( x\right) }{%
\lambda \left( x\right) }df\left( e_{3}\left| _{x}\right. \right) \medskip ,
\\
dN_{2}\left( X_{2}\right) &=&-\frac{u\left( x\right) }{\lambda \left(
x\right) }df\left( e_{1}\left| _{x}\right. \right) -\frac{v\left( x\right) }{%
\lambda \left( x\right) }df\left( e_{3}\left| _{x}\right. \right) .
\end{eqnarray*}
Denote by $\widetilde{A}_{1}$, $\widetilde{A}_{2}$ the shape operators of $%
\widetilde{\xi }$ at $\pi \left( x\right) $ corresponding to the directions $%
N_{1}$ and $N_{2}$. Taking into account the above relations, from Weingarten
formulas it follows that at $\pi \left( x\right) $ we have 
\begin{equation*}
\widetilde{A}_{1}\sim \frac{1}{\lambda \left( x\right) }\left( 
\begin{array}{ll}
1 & \ \ 0 \\ 
0 & -1
\end{array}
\right) ,\quad \widetilde{A}_{2}\sim \frac{1}{\lambda \left( x\right) }%
\left( 
\begin{array}{ll}
-v\left( x\right) & -u\left( x\right) \\ 
-u\left( x\right) & \ \ v\left( x\right)
\end{array}
\right) ,
\end{equation*}
with respect to the orthonormal base $\left( X_{1},X_{2}\right) $. So the
immersion $\widetilde{\xi }:V\rightarrow \mathbb{S}^{4}$ is a minimal
immersion whose Gaussian curvature $\mathcal{K}$ and normal curvature $%
\mathcal{K}_{n}$ are given by 
\begin{equation}
\mathcal{K}=1-\frac{1+u^{2}+v^{2}}{\lambda ^{2}},\ \ \mathcal{K}_{n}^{2}=%
\frac{4u^{2}}{\lambda ^{4}}>0.
\end{equation}
This completes the proof.
\end{proof}

\noindent We shall use in the proof of our\ theorem a result due to Cheng
and Yau \cite{CY} that we recall in the following lemma.

\begin{lemma}
\textit{Let }$M^{n}$\textit{\ be an }$n$\textit{-dimensional, }$n\geq 2$,%
\textit{\ complete Riemannian manifold with Ricci curvature }$Ric\geq
-\left( n-1\right) k^{2}$,\textit{\ where }$k$ \textit{is a positive constant%
}.\textit{\ Suppose that }$h$\textit{\ is a smooth non-negative function on }%
$M^{n}$ \textit{satisfying} 
\begin{equation*}
\Delta h\geq ch^{2},
\end{equation*}
\textit{where }$c$ \textit{is\ a positive constant and }$\Delta $ \textit{%
stands for the Laplacian operator. Then }$h$ \textit{vanishes identically.}
\end{lemma}

\textit{Proof of the\ Theorem. }After passing to the universal covering
space of $M^{3},$ we may suppose that $M^{3}$ is simply connected. Since $%
M^{3}$ is simply connected, it is oriented and the standard monodromy
argument allows us to define a global orthonormal frame field $\left(
e_{1},e_{2},e_{3}\right) $ of principal directions. Our assumptions imply
that $M^{3}$ has three distinct principal curvatures $\lambda >0>-\lambda .$
The functions $u$ and $v$ are well defined on entire $M^{3}$. By virtue of
Lemma 2, we may assume that $u>0$.

Let $V$ be the quotient space of leaves of $e_{2}$. According to Lemma 3 the
immersion $\widetilde{\xi }:V\rightarrow \mathbb{S}^{4}$ is minimal with
normal curvature nowhere zero. Denote by $\mathcal{N}$ the unit normal
bundle of $\widetilde{\xi }$ in $\mathbb{S}^{4}$. Then the polar map $\Psi :%
\mathcal{N\rightarrow }\mathbb{S}^{4}$, $\Psi \left( x,w\right) =w$, is an
immersion. Consider the map $\tau :M^{3}\rightarrow \mathcal{N},$ $\tau
\left( x\right) =\left( \pi \left( x\right) ,f\left( x\right) \right) .$
Because $\Psi \circ \tau =f$ it follows that$\ \tau $ is a local isometry
and $\Psi \left( \mathcal{N}\right) \equiv f\left( M^{3}\right) $. Hence $%
f\left( M^{3}\right) $ is the image of the polar map associated with $%
\widetilde{\xi }$.

Using $\left( 2.3\right) $ and the fact that $u$ and $v$ are harmonic, we
obtain 
\begin{eqnarray*}
\frac{1}{2}\Delta \left( \left( u-1\right) ^{2}+v^{2}\right) &=&\left|
\nabla u\right| ^{2}+\left| \nabla v\right| ^{2} \\
&\geq &\left( e_{2}\left( u\right) \right) ^{2}+\left( e_{2}\left( v\right)
\right) ^{2} \\
&=&4u^{2}v^{2}+\left( v^{2}-u^{2}+1\right) ^{2} \\
&=&\left( \left( u-1\right) ^{2}+v^{2}\right) ^{2}+4u\left( \left(
u-1\right) ^{2}+v^{2}\right) .
\end{eqnarray*}
Therefore we have 
\begin{equation*}
\Delta \left( \left( u-1\right) ^{2}+v^{2}\right) \geq 2\left( \left(
u-1\right) ^{2}+v^{2}\right) ^{2}.
\end{equation*}
In view of our assumptions, the Ricci curvature of $M^{3}$ is bounded from
below. Appealing to Lemma 4, we deduce that $\left( u-1\right) ^{2}+v^{2}$
is identically zero. Consequently $u\equiv 1$ and $v\equiv 0$. From $\left(
3.4\right) $ it follows that $\left( 1-\mathcal{K}\right) ^{2}-\mathcal{K}%
_{n}^{2}=0$ and so\ the isometric immersion $\widetilde{\xi }$ is
superminimal.

Suppose now that $0<\inf S\leq \sup S<\infty $. At first we will show that $V
$ is complete with respect to the metric $\left\langle \ ,\ \right\rangle $
induced by $\widetilde{\xi }$.\ Arguing indirectly, assume that $V$ is not
complete. Then, there exist a divergent curve $a:[0,\infty )\rightarrow V$
with finite length. The curve $c:[0,\infty )\rightarrow \mathcal{N},$ $%
c\left( t\right) :=\left( a\left( t\right) ,\eta \left( t\right) \right) $,
where $\eta \left( t\right) $ is a unit normal vector field parallel, along $%
a$, with respect to the normal connection of $\widetilde{\xi }$, is
divergent. Because the induced metric on $\mathcal{N}$ by $\Psi $ is
complete, $c$\ has infinite length. Moreover, we have 
\begin{eqnarray*}
\Psi \left( c\left( t\right) \right)  &=&\eta \left( t\right) , \\
d\Psi \left( c^{\prime }\left( t\right) \right)  &=&\frac{\overline{\nabla }%
\eta }{dt}=-\widetilde{A}_{\eta }\left( a^{\prime }\left( t\right) \right) ,
\end{eqnarray*}
where $\overline{\nabla }$ is the Levi-Civita connection on $\mathbb{R}^{5}$
and $\widetilde{A}_{\eta }$ is\ the shape operator of $\widetilde{\xi }$
associated with $\eta $. Because $\inf S>0,$ we have 
\begin{equation*}
\left\| \widetilde{A}_{\eta }\right\| \leq \sqrt{2\left( 1-\mathcal{K}%
\right) }=2\sqrt{\frac{2}{S}}\leq 2\sqrt{\frac{2}{\inf S}}.
\end{equation*}
Then, 
\begin{eqnarray*}
\int_{0}^{\infty }\left| c^{\prime }\left( t\right) \right| dt
&=&\int_{0}^{\infty }\left| \widetilde{A}_{\eta }\left( a^{\prime }\left(
t\right) \right) \right| dt \\
&\leq &\int_{0}^{\infty }\left\| \widetilde{A}_{\eta }\right\| \left| \left(
a^{\prime }\left( t\right) \right) \right| dt \\
&\leq &2\sqrt{\frac{2}{\inf S}}\int_{0}^{\infty }\left| a^{\prime }\left(
t\right) \right| dt<\infty \text{,}
\end{eqnarray*}
which leads to a contradiction. Hence $V$ must be complete. Endow now $V$
with the conformal metric 
\begin{equation*}
\widehat{\left\langle \ ,\ \right\rangle }=\left( 1-\mathcal{K}\right) ^{%
\frac{1}{3}}\left\langle \ ,\ \right\rangle .
\end{equation*}
Note that $\widehat{\left\langle \ ,\ \right\rangle }$ is complete, since $1-%
\mathcal{K}\geq \frac{4}{\sup S}$.\ The Gaussian curvature $\widehat{%
\mathcal{K}}$ of the new metric is given by 
\begin{equation*}
\widehat{\mathcal{K}}=\frac{\mathcal{K}}{\left( 1-\mathcal{K}\right) ^{\frac{%
1}{3}}}-\frac{\Delta \log \left( 1-\mathcal{K}\right) ^{\frac{1}{3}}}{%
2\left( 1-\mathcal{K}\right) ^{\frac{1}{3}}},
\end{equation*}
where, here, $\Delta $ stands for the Laplacian with respect to the metric $%
\left\langle \ ,\ \right\rangle $. Inserting in this the identity $\left(
2.1\right) ,$ we find 
\begin{equation*}
\widehat{\mathcal{K}}=\frac{1}{3\left( 1-\mathcal{K}\right) ^{\frac{1}{3}}}%
\geq \frac{1}{3}\left( \frac{\inf S}{4}\right) ^{\frac{1}{3}}.
\end{equation*}
Therefore, $\widehat{\mathcal{K}}$ is bounded away from zero and thus, by
Myers' theorem, $V$ is compact and thus $f\left( M^{3}\right) $ is compact.
In particular, by a result due to Asperti (\cite{A}, Theorem 1), $V$ is
diffeomorphic to $\mathbb{S}^{2}$ or to $\mathbb{RP}^{2}$ and this completes
the proof.$\blacksquare $\medskip 

\noindent \textbf{Remark.} We emphasise that the quotient space $V$ in the
Theorem may be non orientable, although $M^{3}$ is orientable. To illustrate
this, consider the \textit{Veronese surface }$g:\mathbb{S}%
_{1/3}^{2}\rightarrow \mathbb{S}^{4}$, 
\begin{equation*}
g\left( x,y,z\right) =\left( \frac{xy}{\sqrt{3}},\frac{xz}{\sqrt{3}},\frac{yz%
}{\sqrt{3}},\frac{x^{2}-y^{2}}{2\sqrt{3}},\frac{x^{2}+y^{2}-2z^{2}}{6}%
\right) \text{,}
\end{equation*}
which induces an isometric embedding $\widetilde{g}:\mathbb{RP}%
^{2}\rightarrow \mathbb{S}^{4}$. Then, the unit normal bundle $\mathcal{N}$
of $\widetilde{g}$ is compact and the polar map $\Psi :\mathcal{N}\mathbb{%
\rightarrow S}^{4}$ provides a minimal isoparametric hypersurface with
principal curvatures $\sqrt{3}$, $0$, $-\sqrt{3}$, the so called Cartan
hypersurface. Since the immersion $\Psi $ admits $\widetilde{g}$ as global
normal vector field, the manifold $\mathcal{N}$ is orientable.

Concluding, we pose the following question: \textit{Does there exist a
complete minimal hypersurface }$f:M^{3}\rightarrow \mathbb{S}^{4}$\textit{\
with }$K=0$\textit{\ and }$S>0$\textit{\ whose Gauss image is not
superminimal?} Of course, if such an example exists, its $S$ must be
unbounded.

\end{document}